\documentclass[12pt]{article}
\usepackage{amsfonts,amsmath,amssymb}
\usepackage{enumerate}
\usepackage{theorem}
\numberwithin{equation}{section}

\newtheorem{theorem}{Theorem}[section]
\newtheorem{lemma}[theorem]{Lemma}
\newtheorem{proposition}[theorem]{Proposition}

\newtheorem{definition}{Definition}[section]
\newtheorem{corollary}[theorem]{Corollary}

\newcommand{\la}{\lambda}
\newcommand{\cm}{\mathcal{M}}
\newcommand{\mi}{\mu}
\newcommand{\ra}{\rightarrow}
\newcommand{\R}{\mathbb{R}}
\newcommand{\bi}{\beta}

\begin{document}

\title{Extremal problems related to maximal dyadic-like operators}
\author{Eleftherios N.Nikolidakis\\
{\small Department of Mathematics and Statistics},\vspace*{-0.15cm}\\
{\small University of Cyprus 20537, CY 1678, Nicosia, Cyprus}}
\footnotetext{\hspace{-0.5cm}} \footnotetext{\hspace{-0.5cm}E-mail
address: nikolidakis.lefteris@ucy.ac.cy}
\date{}
\maketitle

\begin{abstract}
We obtain sharp estimates for the localized distribution function of the dyadic maximal function $\mathcal{M}_d\phi$.
when $\phi$ belongs to $L^{p,\infty}$. Using this we obtain sharp estimates for the quasi-norm of $\mathcal{M}_d\phi$ in
$L^{p,\infty}$ given the localized $L^1$-norm and certain weak $L^p$-conditions
\medskip



\end{abstract}

\section{Introduction}
The dyadic maximal operator on $\mathbb{R}^n$ is defined by
\begin{equation}
\mathcal{M}_d\phi (x)= \sup\left\{ \frac{1}{|Q|}\int_Q |\phi (u)| du : x\in Q, Q\subseteq \mathbb{R}^n\: \mbox{is a dyadic cube} \right\}
\end{equation}
for every $\phi\in L^1_{loc}(\mathbb{R}^n)$, where the dyadic cubes are those formed by the grids $2^{-N}\mathbb{Z}^n$ for
$N=0,1,2,...$. It is known that it satisfies the following weak type $(1,1)$ inequality
\begin{equation}
|\left\{x\in\mathbb{R}^n : \mathcal{M}_d\phi (x)>\lambda       \right\}|\leq \frac{1}{\lambda}
\int_{ \left\{ \mathcal{M}_d\phi (x)>\lambda \right\}    }|\phi (u) |du
\end{equation}
for every $\phi\in L^{1}(\mathbb{R}^n)$ and every $\lambda >0$ from which we  easily get the following $L^p$ inequality
\begin{equation}
|| \mathcal{M}_d\phi  ||_p\leq \frac{p}{p-1}|| \phi ||_p  \label{best_poss}.
\end{equation}
One way of studying such maximal operators is the computation of the
so-called Bellman functions related to them, which reflect certain
deeper properties of these maximal operators. Such functions related
to inequality (\ref{best_poss}) are precisely evaluated
\cite{Melas_monos}, \cite{Melas_koino}. Actually if we define for
any $p>1$

\begin{equation}
B_p(f,F)=\sup\left\{ \frac{1}{|Q|}\int_Q (\mathcal{M}_d\phi)^p du :\:  \frac{1}{|Q|}\int_Q\phi=f, \:\frac{1}{|Q|}\int_Q\phi^p=F
 \right\} \label{belman}
\end{equation}
where $Q$ is a fixed dyadic cube, $\phi$ is non-negative in $L^p(Q)$ and $F,f$ satisfy $0\leq f \leq F^{1/p}$. It is proved
in \cite{Melas_monos} that $B_p(f,F)=F\omega_p (\frac{f^p}{F})^p$ where $\omega_p :[0,1]\rightarrow [1,\frac{p}{p-1}]$ is the
inverse function of $H_p(z)=-(p-1)z^p+pz^{p-1}$.Using the above result it is possible to compute more complicated functions
related to maximal dyadic-like operators. \newline\newline
The case where $p<1$ is studied in \cite{Melas_koino} where it was completely solved.\newline\newline
One may look (\ref{belman}) as an extremum problem which reflects the deeper structure of dyadic-maximal operators. Certain
other extremum problems arise in this spirit. Some of these are the computation of the following functions

\begin{align}
B(f,F)=\sup\bigg\{|| \mathcal{M}_{\mathcal{T}}\phi||_{p,\infty}  :\:\phi\ge0,  \int_X\phi d\mu=f, \: &|||\phi|||_{p,\infty}=F
 \bigg\},\:\nonumber\\
 &\hspace*{0.8cm} 0<f\leq F  \label{extremum1}
\end{align}
and

\begin{align}
B_1 (f,F)=\sup\bigg\{|| \mathcal{M}_{\mathcal{T}}\phi||_{p,\infty}  :\:\phi\ge0,  \int_X\phi d\mu=f, \: &||\phi||_{p,\infty}=F
 \bigg\}, \nonumber\\
 &\hspace*{-0.4cm} 0<f\leq \frac{p}{p-1}F  \label{extremum2}
\end{align}
where $(X,\mu)$ is a non-atomic probability measure space, $\mathcal{T}$ a tree on $X$ and $\mathcal{M}_{\mathcal{T}}$ the
corresponding maximal operator as it will defined in the sequel. Additionally

\begin{equation}
||\phi||_{p,\infty}=\sup\left\{  \lambda \mu \left(  \left\{\phi  >\lambda \right\}  \right)^{1/p} :\:  \lambda >0 \right\}
\label{eq1.7}
\end{equation}
is the standard quasi-norm in $L^{p,\infty}$ and
\begin{eqnarray}
sup\left\{ \mu \left( E\right)^{-1+1/p}\int_E|\phi|d\mu : E\subseteq X \mbox{ {\small measurable, such that} } \mu(E)>0 \right\}
\label{eq1.8}
\end{eqnarray}
is  an equivallent norm in $L^{p,\infty}$.

In this article we exactly compute these functions. In fact we set.
\begin{eqnarray}
B(f,F,\la)=\sup\{\mu(\{\cm\phi\ge\la\}):\:\phi\ge0, \int_X\phi d\mu=f, \:|||\phi|||_{p,\infty}\le F\}.
\label{eq1.9}
\end{eqnarray}
for $0<f\le F$, $\la>0$ and
\begin{eqnarray}
B_1(h,F,\la)=\sup\left\{\mu\left(\left\{\cm\phi\ge\la\right\}\right):\:\phi\ge0, \int_X\phi d\mu=f,\|\phi\|_{p,\infty}
=F\right\}  \label{eq1.10}
\end{eqnarray}
for $0<f\le\frac{p}{p-1}F$, $\la>0$ and we compute them. After that
the computation of (\ref{extremum1}) and (\ref{extremum2}) is an
easy task.

Additionally it is not difficult to compute
\begin{equation}
B_2(f,F,\la)=\sup\left\{\mu\left(\left\{\cm\phi\ge\la\right\}\right):\:\phi\ge0, \int_X
\phi d\mu=f, |||\phi|||_{p,\infty}=F\right\}  \label{eq1.11}
\end{equation}
for the range $f<\la\le(F^p/f)^{1/(p-1)}$, by the method that
(\ref{eq1.10}) is computed, while the proof remains the same for
$\la>(P^p/f)^{1/(p-1)}$, as it is in Theorem 3.1. The result is
that:
\[
B_2(f,F,\la)=B(f,F,\la).
\]

Related problems are studied in \cite{koino2}. In all these problems
the corresponding functions are independent from the particular tree
$\mathcal{T}$ and the measure space $(X,\mu)$.

\section{Some general facts}

Let $(X,\mu)$ be an non-atomic probability measure space. We begin with the following
\begin{lemma}
Let $ \phi:X \rightarrow \mathbb{R}^+  $ measurable and $I\subseteq
X$ be measurable with $\mu(I)>0$. suppose that $Av_I(\phi)=
\frac{1}{\mu(I)}\int_I\phi d\mu=s$. Then for every $\beta\in (0,
\mu(I)]$ there exists measurable set $E_\beta\subseteq I$ with
$\mu(E_\beta)=\beta$ such that
$\frac{1}{\mu(E_\beta)}\int_{E_\beta}\phi d\mu=s$.
\end{lemma}
{\bf Proof.}\newline
Consider the measure space $(I,\mi/I)$ and define $\psi:I\ra\R^+$
with $\psi(t)=\phi(t)$, $t\in I$, so that $\psi=\phi/I$.

Then, if $\psi^\ast:[0,\mi(I)]\ra\R^+$ is the decreasing
rearrangement of $\psi$, we have that
\begin{eqnarray}
\frac{1}{\bi}\int^\bi_0\psi^\ast(u)du\ge\frac{1}{\mi(F)}\int^{\mi(I)}_0
\psi^\ast(u)du=s\ge\frac{1}{\bi}\int^{\mi(I)}_{\mi(I)-\bi}\psi^\ast(u)du.
\label{eq2.1}
\end{eqnarray}
The inequalities are obvious because $\psi^\ast$ is decreasing while
the equality is true because of
\[
\int^{\mi(I)}_0\psi^\ast(u)du=\int_I\phi d\mi.
\]
From (\ref{eq2.1}) it is deduced that there exists a $r\ge0$ such
that $\bi+r\le\mi(I)$ with
\begin{eqnarray}
\frac{1}{\bi}\int^{\bi+r}_r\psi^\ast(u)du=s.  \label{eq2.2}
\end{eqnarray}
It is now easily seen that there exists an $E_\bi$ measurable subset
of $I$ such that
\begin{eqnarray}
\mi(E_\bi)=\bi \ \ \text{and} \ \ \int_{E_\bi}\phi d\mi=\int^{\bi+r}_r\psi^\ast
(u)du.   \label{eq2.3}
\end{eqnarray}
(\ref{eq2.2}) and (\ref{eq2.3}) now give the conclusion of Lemma
2.1.
\bigskip
\newline
We prove now the following
\begin{proposition}
Let $ \phi:X \rightarrow \mathbb{R}^+ $ measurable, $\int_X\phi
d\mu=f,\: ||| \phi |||_{p,\infty}=F ( f\leq F)$ . Let $I\subseteq X$
be measurable with $\mu(I)>0$ such that $Av_I(\phi)=s$. Define $\psi
: X\rightarrow \mathbb{R}^+ $ with
\begin{equation}
\psi (t)=\left\{
\begin{array}{ll}
\phi (t) , & t\in X\backslash I \\
        s, & t\in I
\end{array}\right.
\end{equation}
Then $\psi\in L^{p,\infty}$, $|||\psi  |||_{p,\infty}\leq F$ and
$\int_X\psi d\mu=f.$
\end{proposition}
{\bf Proof.}\newline It is clear that $\int_X\psi d\mu=\int_X\phi
d\mu=f$.Let now $E\subseteq X$ be measurable with $\mu(E)>0$. We
prove that $\mu \left( E\right)^{-1+1/p}\int_E\psi d\mu\leq F $. We
write $E=E_1\cup E_2$ with $E_1\subseteq X\backslash I,E_2\subseteq
I$ and $\mu(E_1)+\mu(E_2)>0$. Then
\begin{equation}
\mu \left( E\right)^{-1+1/p}\int_E\psi d\mu =\frac{ \int_{E_1}\psi d\mu+\int_{E_2}\psi d\mu }{ ( \mu(E_1)+\mu(E_2) )^{1-1/p} } =
\frac{ \int_{E_1}\phi d\mu+s\mu(E_2) }{ ( \mu(E_1)+\mu(E_2) )^{1-1/p} } \label{fr1}
\end{equation}
Using Lemma 1 we have that $s\leq \frac{1}{\mu(E_\beta)}\int_{E_\beta}\phi d\mu$ for some $E_\beta$ measurable subset of $I$
with $\mu(E_\beta)=\mu(E_2)$. Hence from (\ref{fr1}) we deduce that
\begin{equation}
\begin{array}{rl}
\mu \left( E \right)^{-1+1/p} \int_E \psi d\mu \leq &  \frac{ \int_{ E_1 }\phi d\mu+\int_{ E_\beta }\phi d\mu }{ ( \mu(E_1)+\mu(E_2) )^{1-1/p}   }\\
                                                =   &  \mu \left( E_1\cup E_\beta  \right)^{-1+1/p} \int_{ E_1 \cup E_\beta } \phi d\mu \\
                                              \leq  & ||| \phi |||_{p,\infty}=F
\end{array}
\end{equation}
so that $||| \psi |||_{p,\infty}\leq F$ and the proposition is proved.\newline
\newline
Let now $(X,\mu)$ be a non-atomic probability measure space. Two measurable subsets $A,B$ of $X$ will be call {\it almost disjoint}
if $\mu(A\cap B)=0$. We  give now the following
\begin{definition}
A set $ \mathcal{T} $of measurable subsets of $X$ will be called a tree if the following conditions are satisfied:\newline
i) $X\in\mathcal{T}$ and for every $I\in \mathcal{T}$ we have $\mu (I)>0$ \newline
ii) For every $I\in \mathcal{T} $ there corresponds a finite subset $\mathcal{C}(I)\subseteq \mathcal{T}$ containing at least two elements
such that:

(a)the elements of $\mathcal{C}(I)$ are pairwise almost disjoint subsets of $I$.

(b)$I= \cup \mathcal{C}(I)$\newline
iii) $\mathcal{T}=\cup_{m\geq (0)} \mathcal{T}_{(m)} $ where $\mathcal{T}_{0}=\{ X \} $ and $\mathcal{T}_{(m+1)}=\cup_{I\in \mathcal{T}_{(m)}} \mathcal{C}(I) $ \newline
iv) $\lim_{m\rightarrow \infty}\mbox{sup}_{ I\in \mathcal{T}_{(m)} }\mu (I)=0$
\end{definition}
From \cite{Melas_monos} we have the following
\begin{lemma}
For every $I\in\mathcal{T}$ and every $a$ such that $0<a<1$ there exists a subfamily $\mathcal{F}(I) \subseteq \mathcal{T}$ consisting of
pairwise almost disjoint subsets of $I$ such that
\begin{equation*}
\mu \left( \cup_ { J\in\mathcal{F}(I) } J \right)=\sum_{J\in\mathcal{F}(I)    }\mu (J)=(1-a)\mu(I).
\end{equation*}
\end{lemma}
Let now $(X,\mu)$ be a non-atomic probability measure space and $\mathcal{T}$ a tree as in definition (1.1). We define
the associated
maximal operator to the tree $\mathcal{T}$ as follows: Let $\phi$ be
 function on $X$ such that $\phi\in L^1(X,\mu)$ then
\begin{equation}
\mathcal{M}_{\mathcal{T}}\phi (x)=sup\left\{ \frac{1}{\mu (I) }\int_I|\phi |d\mu : x\in I\subseteq \mathcal{T} \right\}
\end{equation}
for every $x\in X$. Due to proposition (1.2) in order to find
\begin{equation*}
sup\left\{ \mu(\left\{ \mathcal{M}_{ \mathcal{T} }\phi\geq \lambda \right\}): \phi\geq 0,\int_X\phi d\mu=f,|||\phi |||_{p,\infty}\leq F \right\}
\end{equation*}
we may assume that $\phi=\lambda$ on the set where $\mathcal{M}_{
\mathcal{T} }\phi\geq \lambda $.

As for the domain of the respective extremal problem it is easy to
see that it is as mentioned in equations (\ref{extremum1}) and
(\ref{extremum2}). More precisely there exists a non negative $\phi$
and not equal to the zero function such that $\int_X\phi d\mi=f$ and
$|||\phi|||_{p,\infty}=F$ if and only if $0<f\le F$ for
(\ref{extremum1}) and analogously for (\ref{extremum2}).
\section{The first extremal problem}

\begin{theorem}
Let $(X,\mu)$ be a non atomic probability measure space, $\mathcal{T}$ a tree on the measure space $X$ and $\mathcal{M}_{\mathcal{T}}=\mathcal{M}$
the associated maximal operator, then the following holds
\begin{equation*}
\begin{array}{rl}
B(f,F,\lambda)= & sup\left\{  \mu\left(\left\{\mathcal{M}\phi\geq \lambda  \right\} \right): \phi\geq 0,\int_X\phi d\mu=f,|||\phi |||_{p,\infty}\leq  F   \right\} \\
              = & min(1,f/\lambda, F^p/\lambda^p)
\end{array}
\end{equation*}
where $0<f\leq F,\lambda >0$.
\end{theorem}
{\bf Proof.} \newline
 We  calculate
\begin{equation*}
min(1,f/\lambda, F^p/\lambda^p)=\left\{
\begin{array}{ll}
1,& 0<\lambda\leq f \\
f/\lambda , & f<\lambda\leq  \\
F^p/\lambda^p, & \left(F^p/f \right)^{1/(p-1)} < \lambda
\end{array}
\right.
\end{equation*}
For $0<\lambda\leq f$ it is implied $ \mu( \left\{ \mathcal{M}\phi
\geq \lambda \right\} )=\mu(X)=1$ for any $\phi$ such that
$\int_X\phi d\mu =f$. For $f< \lambda \leq \left(F^p/f
\right)^{1/(p-1)} $ it is proved in \cite{koino2} that
\begin{equation}
sup\left\{  \mu\left(\left\{\mathcal{M}\phi\geq \lambda  \right\} \right): \phi\geq 0,\: \int_X\phi d\mu=f,\: \int_X\phi^p d\mu = F^p   \right\}=f/\lambda
\end{equation}
since $L^p(X,\mu)\subset L^{p,\infty}(X,\mu)$ and $\int_X\phi^pd\mu =F^p$ implies $|||\phi |||_{p,\infty}\leq F$ we lead to $B(f,F,\lambda)= f/\lambda$ due to
\begin{equation*}
 \mu( \left\{ \mathcal{M\phi} \geq \lambda \right\} )\leq \frac{\int_{ \left\{ \mathcal{M\phi} \geq \lambda \right\} }\phi} {\lambda}\leq f/\lambda, \:\forall \lambda >0
\end{equation*}
We prove now the rest of the theorem by showing that $B(f,F,\lambda)= F^p/\lambda^p$, when $\lambda>\left(F^p/f \right)^{1/(p-1)} $. Without loss of
generality we suppose that $F=1, f\leq 1$. We first prove that $B(f,F,\lambda)\leq  1/\lambda^p$ for every $\lambda >0$. For this let $\phi:X\rightarrow \mathbb{R}^+$
 with $|||\phi |||_{p,\infty}\leq 1$, $\int_X\phi d\mu=f$ and $\lambda >0$. If $\lambda<f\leq 1\Rightarrow B(f,F,\lambda)=1<\frac{1}{\lambda^p}$. Let now $\lambda\geq  f$;
 if we set $K=\left\{ \mathcal{M}\phi\geq \lambda \right\}$ we have that $K=\cup_j I_j $ where $I_j \in \mathcal{T}$ almost dijoint
and $\frac{1}{\mu(I_j)}\int_{I_j}\phi d\mu \geq \lambda$ so that $\int_{I_j}\phi d\mu \geq \lambda\mu(I_j), \forall j $. Summing
 over the indices  we obtain
\begin{eqnarray}
\int_{K }\phi d\mu= \int_{\cup_j I_j }\phi d\mu\geq \lambda \mu\left( \cup_j I_j \right) =\lambda \mu(K) \label{phiunion}.
\end{eqnarray}
However, $|||\phi |||\leq 1\Rightarrow \mu(K)^{-1+1/p}\int_K \phi d\mu\leq 1$ and since (\ref{phiunion}) is valid we end up
to $\mu(\left\{\mathcal{M}\phi\geq \lambda \right\})=\mu(K)\leq 1/\lambda^p$. Therefore $B(f,F,\lambda)\leq 1/\lambda^p$ for
every $\lambda >0$. We prove the inverse inequality for $\lambda>\left(1/f\right)^{1/(p-1)}, F=1$. We show it first for the
case where $(X,\mu)=([0,1], |\cdot |)$ where $|\cdot |$ denotes the Lebesque measure. Let $\phi:[0,1]\rightarrow \mathbb{R}^+$
be the function defined by
\begin{equation}
\phi(t)=\left\{
\begin{array}{ll}
(1-1/p)(t+1/\lambda^p)^{-1/p}, & 0\leq t<A \\
0, & A\leq t\leq 1-1/\lambda^p \\
\lambda,& 1-1/\lambda^p<t\leq 1
\end{array}\right.
\end{equation}
where $A=f^{p/(p-1)}-1/\lambda^p\leq 1-1/\lambda^p$. We prove that
$\int_X\phi d\mu=f$, $|||\phi |||_{p,\infty}=1$. For $t\in [0,A]$ we
have that
\begin{equation}
\int_0^t\phi(s)ds=\left(t+1/\lambda^p\right)^{1-1/p}-1/\lambda^{p-1}
\end{equation}
For $A\leq t\leq 1-\frac{1}{\lambda^p}:
\int_0^t\phi(s)ds=f-\frac{1}{\lambda^{p-1}} $ thus $\int_0^1\phi
(s)ds=f$. It is quite clear $|\left\{\phi\geq \lambda \right\}|\geq
|[1-1/\lambda^p,1]|=1/\lambda^p$ and for every $t\in
[0,1-1/\lambda^p]$ we derive that
\begin{equation}
\int_0^t\phi(s)ds\leq \left(t+1/\lambda^p\right)^{1-1/p}-(1/\lambda^p)^{1-1/p} \label{phiunion2}
\end{equation}
therefore for every $E\subset [0,1-1/\lambda^p]$ with $|E|>0 $ we find
\begin{equation}
\int_E\phi(s)ds\leq \int_0^{|E|}\phi(s)ds\leq \left(|E|+1/\lambda^p\right)^{1-1/p}-\left(1/\lambda^p\right)^{1-1/p} \label{phiunion3}
\end{equation}
since  $\phi$ is decreasing on $[0,1-1/\lambda^p]$. \newline\newline
We now prove that $|||\phi |||_{p,\infty}\leq 1$. It is enough to prove that
\begin{equation}
\frac{\lambda |E_1|+ \int_{E_2}\phi (s)ds }{ (|E_1|+|E_2|)^{1-1/p} }\leq 1 \label{lambdaint}
\end{equation}
for every $E_1\subseteq [1-1/\lambda^p,1]$ and every $E_2\subseteq [0,1-1/\lambda^p]$ such that $|E_1|+|E_2|>0$. Equation (\ref{lambdaint}) is equivalent to
\begin{equation}
\begin{array}{rl}
\int_{E_2}\phi (s)ds \leq & \inf\left\{(|E_1|+|E_2|)^{1-1/p}-\lambda |E_1|: |E_1|\leq 1/\lambda^p     \right\} \\
                       =  & \inf\left\{ (a+|E_2|)^{1-1/p}-\lambda a: 0 \leq a\leq 1/\lambda^p     \right\}.
\end{array}
\end{equation}
We define $g:[0,1/\lambda^p]\rightarrow \mathbb{R}^+$ by $g(a)=(a+|E_2|)^{1-1/p}-\lambda a$. Since $g$ is a concave function
of $a$ and $(x+y)^q\leq x^q+y^q$ for
$x,y\geq 0, q<1$, we derive
\begin{equation*}
\begin{array}{rl}
\inf\left\{g(a): 0\leq a \leq 1/\lambda^p     \right\}= & \min\left\{ g(0),g(1/\lambda^p)   \right\} \\
                                                      = & \min\left\{ |E_2|^{1-1/p}, \left(1/\lambda^p+|E_2|     \right)^{1-1/p}-\left(1/\lambda^p \right)^{1-1/p}      \right\} \\
                                                      = & \left(1/\lambda^p+|E_2|    \right)^{1-1/p} -1/\lambda^{p-1}.
\end{array}
\end{equation*}
Therefore in order to prove that $||| \phi|||_{p,\infty}\leq 1$ it is enough to show
\begin{equation*}
\int_{E_2} \phi(s)ds\leq \left(1/\lambda^p+|E_2|    \right)^{1-1/p} -1/\lambda^{p-1}
\end{equation*}
for every $E_2\subseteq [0,1-1/\lambda^p]$ which is true by
(\ref{phiunion3}). Thus we obtain $|||\phi|||_{p,\infty}\leq 1$. In
addition $|||\psi|||_{p,\infty}\geq ||\psi||_{p,\infty}$ for every
$\psi\ \in  L^{p,\infty}$ and $$| \left\{ \phi > \lambda-\epsilon
\right\} |\geq 1/\lambda^p$$ for every $\epsilon >0$ small enough
for  our  specified $\phi$. Therefore $||\phi||_{p,\infty}\geq 1$
consequently $$|||\phi |||_{p,\infty}\geq 1.$$ We have proved that
$||| \phi  |||_{p,\infty}=1$.\newline  We prove now that for
$\lambda > \left( \frac{1}{f}   \right)^{ \frac{1}{p-1} }, F=1$ we
have $B(f,F,\lambda)\geq 1/\lambda^p$ for the general case. Let
$(X,\mu)$ be a non atomic probability space and $\mathcal{T}$ be a
tree on $X$. From Lemma 2.3 there exists a sequence $\left\{ I_j
\right\}_j$ of almost disjoint subsets of $X$ such that $\mu(\cup
I_j)=\sum_j\mu(I_j)=1/\lambda^p$. We set $Z=X\backslash \left(\cup_j
I_j  \right)$. Then $\mu(Z)=1-1/\lambda^p$. Let now $\phi$ be the
function defined  by the following way:  $\phi : [0,1-1/\lambda^p]
\rightarrow \mathbb{R}^+$ with
\begin{equation}
\phi (t)=\left\{ \begin{array}{rl}
                  (1-\frac{1}{p})(t+1/\lambda^p)^{-1/p}, & 0<t\leq A \\
                      0,  & A<t\leq 1-1/\lambda^p
                  \end{array}
                   \right.
\end{equation}
where $A$ is such that $\int_0^A\phi (s)ds=f-1/\lambda^{p-1}$.\newline
Since $\int_0^{1-1/\lambda^p}(1-\frac{1}{p})(s+1/\lambda^p)^{-1/p}ds=1-1/\lambda^{p-1}\geq f-1/\lambda^{p-1}$ this choice of $A$ is possible. \newline
The measure space $(Z,\mu |Z)$ is  non-atomic finite thus there exists $\psi_1 : Z \rightarrow \mathbb{R}^+$ with
$\psi^{*}_1=\phi$, where $\psi^*_1$ is the decreasing rearrangement of $\psi_1$. Define now
$\psi : X\rightarrow \mathbb{R}^+$ by
\begin{equation}
\psi (t)=\left\{
\begin{array}{rl}
\psi_1(t), &  t\in Z         \\
\lambda , &   t\in \cup_j I_j
\end{array}\right.
\end{equation}
Then $\psi$ is a well defined measurable function on $X$ and
\begin{equation}
\begin{array}{rl}
\int_X\psi d\mu = &\int_Z\psi_{1} d\mu +\lambda\mu( \cup_j I_j ) \\
               \: & \: \\
                = & \int_0^{1-1/\lambda^p}\phi (u) du+1/\lambda^{p-1}=f
\end{array}
\end{equation}
In addition for every $E\subseteq Z$
\begin{equation}
\int_E\psi d\mu=\int_E\psi_1 d\mu\leq \int_0^{|E|}\psi^*_1(s)ds=\int_0^{|E|}\phi (u)du\leq (|E|+1/\lambda^p)^{1-1/p}
-1/\lambda^{p-1}
\end{equation}
so as before we can conclude that $|||\psi|||_{p,\infty}\leq 1$. Due
to $|\left\{ \psi=\lambda \right\}|\geq |\cup_j I_j |\geq
1/\lambda^p$ we can see that $||\psi||_{p,\infty}\geq 1$ and from
the inequality $||\psi ||_{p,\infty}\leq ||| \psi |||_{p,\infty}$ we
deduce $|||\psi|||_{p,\infty}=1$. Finally $\mu (\left\{
\mathcal{M}\psi\geq \lambda    \right\}  )\geq 1/\lambda^p$ and this
implies $ B(f,1,\lambda)\geq 1/\lambda^p $ hence $B(f,1,\lambda)=
1/\lambda^p $ for $\lambda > (\frac{1}{f})^{ \frac{1}{p-1} }$.
Theorem 4.1 is proved.\newline\newline A direct consequence of the
Theorem 4.1 is the following

\begin{corollary}
Let $0<f\leq F$ and
$$B(f,F)=\sup\left\{||\mathcal{M}\phi||_{p,\infty}:\int_X\phi d\mi=f,
|||\phi |||_{p,\infty}=F \right\}.$$ Then $B(f,F)=F$
\end{corollary}
{\bf Proof.} \newline It is easy to see that in the definition of
$B(f,F,\lambda)$ we can replace $\mu (\left\{\mathcal{M}\phi \geq
\lambda \right\})$ by $ \mu (\left\{\mathcal{M}\phi >  \lambda
\right\}) $. By theorem (1.4) we derive $\mu (\left\{\mathcal{M}\phi
>  \lambda \right\})\leq \frac{F^p}{ \lambda^p} $ for every $\phi$
such that $|||\phi |||_{p,\infty}=F$ and for every $\lambda >0$;
therefore $\lambda^p \mu \left(  \left\{ \mathcal{M}\phi >\lambda
\right\}   \right)\leq F^p$, which means $||
\mathcal{M}\phi||_{p,\infty}\leq F$ if $|||\phi |||_{p,\infty}=F$.
However, for $\lambda >(F^p/f)^{1/(p-1)}$ it is implied
$B(f,F,\lambda)=F^p/\lambda^p$ and the supremum is attained for some
$\phi $ such that $|||\phi |||_{p,\infty}=F$ and $\int_X\phi d\mi=f$
hence $B(f,F)=F$ for $0<f\leq F$.

\section{The second extremal problem}

We now prove the following:
\begin{theorem}
Let $(X,\mu )$ be a non atomic probability measure space,
$\mathcal{T}$ a tree on $X$, $\mathcal{M}_{ \mathcal{T} }=
\mathcal{M}$ the associated maximal operator. Then for $0<f\leq
\frac{p}{p-1}\,F$ we have that $$B_1(f,F,\lambda )=\sup \left\{\mu
(\ \{\mathcal{M}\phi\geq \lambda \}):\int_X\phi \,d\mu =f,\: ||\phi
||_{p,\infty}=F\right\}$$
\begin{equation*}
=\min(1,f/\lambda,k^p F^p/\lambda^p)=\left\{
\begin{array}{ll}
1, & 0<\lambda<f \\
f/\lambda, & f\leq \lambda < (F^pk^p/f)^{1/(p-1)} \\
k^pF^p/\lambda^p, & \lambda\geq (F^pk^p/f)^{1/(p-1)}
\end{array}\right.
\end{equation*}
where $k=p/(p-1)$ and $||\phi||_{p,\infty}=\sup \left\{ \lambda \mu\left( \left\{ \phi>\lambda\right\}   \right)^{1/p}: \: \lambda >0 \right\}     $
\end{theorem}
{\bf Proof.} \newline We assume that $F=1$, $f\leq k$ and we search
for $B_1(f,1,\lambda)$. We observe that $||\phi||_{p,\infty}=1$
implies $|||\phi|||_{p,\infty}\leq k$ since for  every  $\psi\in
L^{p,\infty}$  it  is  true  that $|||\psi|||_{p,\infty}\leq
\frac{p}{p-1}||\psi||_{p,\infty}$ Therefore
\begin{equation}
\mu\left( \left\{ \mathcal{M}\phi\geq \lambda \right\} \right)^{-1+1/p} \int_{ \left\{ \mathcal{M}\phi\geq \lambda \right\} }\phi d\mu\leq k \label{intnormp}
\end{equation}
however, $\int_{ \left\{ \mathcal{M}\phi\geq \lambda \right\} }\phi d\mu \geq \lambda \cdot \mu\left( \left\{ \mathcal{M}\phi\geq \lambda \right\}\right)$
thus (\ref{intnormp}) yields $\mu\left(\left\{\mathcal{M}\phi\geq \lambda  \right\}   \right)\leq k^p/\lambda^p$ for every
$\lambda >0$. Hence for $\lambda > \left( k^p/f \right)^{1/(p-1)}$,
\begin{equation}
B_1(f,1,\lambda)\leq k^p/\lambda^p \label{normB2}
\end{equation}
We proceed  proving  the reverse inequality of (\ref{normB2}); Lemma
2.3 implies the existense of a sequence $I_j$ of pairwise almost
disjoint subsets of $X$ such that $I_j\in \mathcal{T}$ for every $j$
and $\sum_j \mu(I_j)=\mu \left( \cup_j I_j  \right)=k^p/\lambda ^p$.
For all values of $\lambda$ satisfying $k^p/\lambda^p<f/\lambda<1$
the particular choice of $(I_j)_j$ is possible. Since
\begin{equation}
\int_0^{k^p/\lambda^p}t^{-1/p}dt=k^p/\lambda^{p-1}=\lambda\sum_j\mu(I_j)
\end{equation}
it is not difficult to see using Lemma 2.1 repeatedly that there
exists a sequence $(A_j)_j$ of measurable pairwise disjoint subsets
of $[0,k^p/\lambda^p]$ such that $\sum_j |A_j|=k^p/\lambda^p$ and
$\int_{A_j}t^{-1/p}dt=\lambda |A_j|$ with $|A_j|=\mu (I_j)$ for any
$j$. For every $j$ we set $\xi_j:A_j\rightarrow \mathbb{R}^+$
defined as $\xi_j(t)=t^{-1/p}$ for $t\in A_j$. Let
$\psi_j:[0,|A_j|]\rightarrow \mathbb{R}^+$ with $\psi_j=\xi^*_j$,
the decreasing rearrangement of $\xi_j$ and since $(I_j, \mu|I_j)$
is a non-atomic finite measure space, there exists measurable
$\phi_j :I_j\rightarrow \mathbb{R}^+$ such that $\phi_j^*=\psi_j$.
Therefore
\begin{equation*}
\int
_{I_j}\phi_j d\mu=\int_0^{|A_j|}\psi_j(t) dt =\int_{A_j}\xi_j(t) dt =\int_{A_j}t^{-1/p}dt=\lambda\mu(I_j)
\end{equation*}
thus
\begin{equation*}
Av_{I_j}(\phi)=\frac{1}{\mu (I_j)}\int_{I_j}\phi d\mu=\lambda.
\end{equation*}
Since $\left\{I_j\right\}_j$ is almost disjoint there exists a
well-defined function $\phi^\prime: \cup_j I_j \rightarrow
\mathbb{R}^+ $ such that $\phi^\prime=\phi_j$ $\mu$-almost
everywhere on $I_j$. Then for every $\theta >0$ with $1/\theta^p\leq
k^p/\lambda^p$ we deduce that
\begin{equation*}
\begin{array}{rl}
\mu(\left\{ \phi^\prime>\theta\right\})=& \sum_j\mu(\left\{ \phi^\prime>\theta \right\}\cap I_j)= \sum_j|\left\{\phi^*_j>\theta \right\}| \\
                                       =&  \sum_j|\left\{ \psi_j>\theta \right\}| =|\left\{ t>0: t^{-1/p}>\theta \right\}| \\
                                       =& 1/\theta^p
\end{array}
\end{equation*}
For $\theta>0$ such that $1/\theta^p>k^p/\lambda^p$,
\begin{equation*}
\mu(\left\{\phi^\prime>\theta \right\})\leq \mu ( \cup_j I_j )=k^p/\lambda^p.
\end{equation*}
The following  holds
\begin{equation*}
\int_{\cup_j I_j }\phi'd\mu=\lambda\sum_j\mu (I_j)=k^p/\lambda^{p-1}
\end{equation*}
We extend the function $\phi'$ to $\phi:X\rightarrow \mathbb{R}^+$ such that $\int_X\phi d\mu=f$,
$||\phi||_{p,\infty}=1$ and $\left\{ \mathcal{M}\phi \geq \lambda  \right\}\supseteq \cup_j I_j$. We proceed
as follows to define $\phi _ \lambda:[0,1-k^p/\lambda^p ]\rightarrow \mathbb{R}^+$ as
\begin{equation*}
\phi_\lambda (t)=\left\{
\begin{array}{ll}
(t+k^p/\lambda^p)^{-1/p} &  ,0\leq t\leq A_1 \\
                       0 & ,A_1<t\leq 1-k^p/\lambda^p
\end{array}
\right.
\end{equation*}
where $A_1$ is such that
\begin{equation}
\int_0^{A_1} \phi_\lambda  (t)dt=f-k^p/ \lambda^{p-1} \label{integral5}
\end{equation}
Since
$\int_0^{1-k^p/\lambda^p}(t+k^p/\lambda^p)^{-1/p}dt=k-k^p/\lambda^{p-1}\geq
f-k^p/\lambda^{p-1}$  , (\ref{integral5}) is possible for some
$A_1\leq 1-k^p/\lambda^p$. Then for
$\theta<(k^p/\lambda^p)^{-1/p}=\la/k$
\begin{equation*}
|\left\{  \phi_\lambda\ge \theta \right\}|\le|\left\{ t\in [0,1-k^p/\lambda^p]:(t+k^p/\la^p)^{-1/p}>\theta \right\}|=\theta^{-p}-k^p/\lambda^p
\end{equation*}
However
\begin{equation*}
\phi_\lambda(t)\leq \lambda/k
\end{equation*}
for every $t\in [0,1-k^p/\lambda^p]$, which gives
\begin{equation*}
\mid\left\{\phi_\lambda >\theta \right\}\mid =0
\end{equation*}
for $\theta>\lambda/k$. Obviously
\begin{equation*}
\mu(X\backslash (\cup_j I_j))=\mu (Z)=1-k^p/\lambda^p
\end{equation*}
so there exists $\phi^{\prime\prime}:Z\rightarrow \mathbb{R}^+$ with$(\phi^{\prime\prime})^*=\phi_\lambda $.
Then
\begin{equation*}
\int_Z\phi^{\prime\prime}d\mu =f-k^p/\lambda^{p-1}
\end{equation*}
and
\begin{equation*}
\mu\left( \left\{ \phi^{\prime\prime}\ge\theta   \right\} \right)\le\left\{
\begin{array}{lr}
0, & \mbox{if}\:\: \theta >\lambda /k \\
\theta^{-p}-k^p/\lambda^p,  &  \mbox{if}\:\: \theta \leq \lambda /k
\end{array}
\right.
\end{equation*}
hence by defining
\begin{equation*}
\phi(t)=\left\{
\begin{array}{ll}
\phi'(t), &  \mbox{if}    \:\: t\in \cup_jI_j    \\
\phi''(t),&  \mbox{if}    \:\: t\in Z
\end{array}
\right.
\end{equation*}
we obtain
\begin{eqnarray}
\int_X \phi d\mu =f,\:\: \parallel \phi \parallel _{p,\infty}=1 \\
\mu (\left\{\mathcal{M}\phi \geq \lambda \right\})\geq k^p/\lambda^p
\end{eqnarray}
therefore  $B_1(f,1,\lambda )=k^p/\lambda^p$ for $\lambda
>(k^p/f)^{1/(p-1)}$. In order to prove the theorem  it is
enough to show that $B_1(f,1,\lambda )=f/\lambda $ for $f<\lambda
\leq (k^p/f)^{1/(p-1)}$ We show it using the following arguments;
for every $\phi \in L^{p,\infty}$ such that $\int_X\phi d\mu =f$, it
is true that $\mu (\left\{\mathcal{M}\phi \geq \lambda \right\})\leq
f/\lambda $, so that $B_1(f,1,\lambda )\leq f/\lambda $. We also
need to prove the reverse inequality;it is easy to see that there
exists $G:[0,f/\lambda ]\rightarrow \mathbb{R}^+$ strictly
decreasing, continous, convex function such that
\begin{equation}
\begin{array}{cc}
\int_0^{f/\lambda}G(t)dt=f,\: G(t)\leq t^{-1/p } \:\:\mbox{for every}\: t \\
\: \\
 \lim_{ t\rightarrow 0^+ } G(t)=+\infty  ,\: \lim_{t\rightarrow +\infty }\left( \mid \left\{ G>t \right\}  \mid t^p\right)=1
\end{array}\label{limscond}
\end{equation}
The existence of such a function $G$ is possible since the following holds
$$\int_0^{f/\lambda }t^{-1/p}dt=\kappa(f/\lambda)^{1-1/p}>f.$$
By Lemma 2.3 there exists a sequence $\left\{I_j\right\}$ of
pairwise almost disjoint elements of $\mathcal{T}$ such that
$\sum_j\mu(I_j)=f/\lambda$. Then, $\int_0^{f/\lambda}G(t)dt=\lambda
\sum_j\mu (I_j)$. Using again Lemma 2.1 repeatedly we see that there
exists a measurable partition of $[0,f/\lambda ]$,
$\left\{A_j\right\}$ such that $\mid A_j\mid =\mu (I_j)$ and
$\int_{A_j}G(t)dt=\lambda \mid A_j\mid $, for every j. We define
$\phi_j:I_j\rightarrow \mathbb{R}^+$ by $\phi_j^*=(G|A_j)^*$. Then
$\int_{I_j}\phi_jd\mu =\int_{A_j}G(t)dt=\lambda \mu (I_j)$, implying
$Av_{I_j}(\phi) =\lambda $ where $\phi =\phi _j$ almost everywhere
on $I_j$ and $\phi =0$ on the complement of $\cup_j I_j$. Then for
every $\theta $ such that $\theta^{-p}\leq f/\lambda $

$$\mu (\left\{\phi >\theta \right\})=\sum_j\mu(\left\{\phi_j>\theta \right\}\cap I_j)=|\left\{G>\theta \right\}|\leq \theta^{-p}$$
However
$$ \mu (\left\{\phi >\theta \right\})\leq \mu (\cup_j I_j)=f/\lambda $$
for every $\theta $ such that $f/\lambda \leq \theta^{-p}$. Thus
$$ \int_X\phi d\mu =f \:\mbox{and}\: ||\phi||_{p,\infty}\leq 1. $$
Due to (\ref{limscond}) $|| \phi||_{p,\infty}=1$. However,
$$\mu (\left\{\mathcal{M}\phi \geq \lambda \right\})\geq \mu (\cup_j I_j)=f/\lambda $$
As a result $B_1(f,1,\lambda)=f/\lambda $ for $f<\lambda \leq
(k^p/f)^{1/(p-1)}$. The theorem is proved.\newline\newline We state
the following
\begin{corollary}
Let $B_1(f,F)=\newline =\sup\left\{||\mathcal{M}\phi
||_{p,\infty}:\phi \geq 0,\int_X\phi d\mu =f,||\phi
||_{p,\infty}=F\right\}$ where $f\leq k F$ and $k =p/(p-1)$. Then
$B_1(f,F)=k F$
\end{corollary}
{\bf Proof.}\newline For any $\phi \geq 0$ such that $\int_X\phi
d\mu =f$ and $||\phi ||_{p,\infty}=F$, we deduce that $|||\phi
|||_{p,\infty}\leq kF$ thus from Corrolary 1.5 $B_1(f,F)\leq k F$.
Theorem 1.4 implies that there exists $\phi \geq 0$ such that
$||\phi ||_{p,\infty}=F$,$\int_X\phi d\mu =f$ and
$\sup\left\{\lambda \mu (\left\{\mathcal{M}\phi \geq \lambda
\right\})^{1/p}:\lambda >0\right\}=k F$ hence $||\mathcal{M}\phi
||_{p,\infty}= k F$. Therefore $B_1(f,F)= k F$ and the corrolary is
proved.

\end{document}